\documentclass{amsart}

\theoremstyle{plain}
\newtheorem{lemma}{Lemma}[section]
\newtheorem{theorem}[lemma]{Theorem}

\newtheorem{proposition}[lemma]{Proposition}

\newtheorem{conjecture}[lemma]{Conjecture}
\newtheorem*{mainTheorem}{Theorem \ref{Main Theorem}}

\theoremstyle{definition}
\newtheorem{definition}[lemma]{Definition}

\numberwithin{equation}{section}

\newcommand{\B}[1]{\mathbf{#1}}
\newcommand{\R}{\mathbb{R}}
\newcommand{\p}{\mathfrak p}
\newcommand{\h}{\mathfrak h}
\newcommand{\g}{\mathfrak g}
\DeclareMathOperator{\diag}{\mathrm{diag}}
\DeclareMathOperator{\SO}{SO}
\DeclareMathOperator{\U}{U}
\DeclareMathOperator{\trace}{trace}
\DeclareMathOperator{\So}{\mathfrak {so}}
\DeclareMathOperator{\spn}{span}

\begin{document}

\title[Nonnegatively curved metrics on Lie groups]{Techniques for
classifying nonnegatively curved left-invariant metrics on compact
Lie groups}

\author{Jack Huizenga}

\address{Department of Mathematics\\ University of Chicago\\
Chicago, IL 60637} \email{huizenga@uchicago.edu}
\thanks{Supported in part by NSF grant DMS-0353634}


\date{August 14, 2006}

\begin{abstract}
We provide techniques for studying the nonnegatively curved
left-invariant metrics on a compact Lie group.  For ``straight''
 paths of left-invariant metrics starting at bi-invariant
metrics and ending at nonnegatively curved metrics, we deduce a
nonnegativity property of the initial derivative of curvature. We
apply this result to obtain a partial classification of the
nonnegatively curved left-invariant metrics on $\SO(4)$.
\end{abstract}

\maketitle

\section{Introduction}
What are the nonnegatively curved left-invariant metrics on a
compact Lie group?  All such metrics on $\SO(3)$ and $\U(2)$ were
classified in \cite{Brown}.  These classifications make use of
various techniques that only work in low dimensions.  For
higher-dimensional groups, the situation becomes complicated, and it
is evident that more tools are necessary to approach the problem
effectively.

In \cite{Tapp}, Tapp proposed to classify the nonnegatively curved
left-invariant metrics on $\SO(4)$ by using infinitesimal
techniques. He considered a path of left-invariant metrics starting
at a bi-invariant metric, and aimed to classify the possibilities
for the initial derivative such that the path looks nonnegatively
curved near the bi-invariant metric. In this article, we provide a
method of passing from this local derivative information to global
information about the entire family of left-invariant metrics.  As a
consequence, we are able to transform Tapp's restrictions on the
derivative of the path into necessary conditions for a
left-invariant metric to be nonnegatively curved.

Let $G$ be a compact Lie group, with bi-invariant metric $h_0$ and
Lie algebra $\g$.  Any left-invariant metric $h$ on $G$ is
determined entirely by its restriction to $\g$, and there is always
some $h_0$-self-adjoint positive definite $\Phi\colon \g\to \g$ such
that
$$h(X,Y) = h_0(\Phi X,Y)$$ holds for $X,Y\in \g$.  In this manner a
smoothly varying family $h_t$ of left-invariant metrics on $G$ can
be identified with a smoothly varying family $\Phi_t$ of
endomorphisms of $\g$.  In this article we are particularly
interested in \emph{inverse-linear} variations of the bi-invariant
metric $h_0$.  We say that $\Phi_t$ is inverse-linear if it is of
the form $\Phi_t= (I-t\Psi)^{-1}$ for some $\Psi\colon \g \to \g$;
of course $\Psi$ is then self-adjoint.  Notice that $\Psi$ is the
derivative of $\Phi_t$ at $t=0$.

When $\Phi_t$ is an inverse-linear variation and $X,$ $Y\in \g$ have
been fixed, we define $\kappa^\Psi(t)$ to be the unnormalized
sectional curvature of $\Phi_t^{-1} X$ and $\Phi_t^{-1}Y$ with
respect to the metric $h_t$.  We omit the superscript $\Psi$ when no
confusion will arise.  For fixed $t_0$, $h_{t_0}$ is nonnegatively
curved if and only if $\kappa(t_0)\geq 0$ for all $X,$ $Y\in \g$.

\begin{definition} The inverse-linear path $\Phi_t$ (or the endomorphism
$\Psi$) is \emph{infinitesimally nonnegative} if for any $X,$ $Y\in
\g$ there is some $\varepsilon>0$ such that the function $\kappa$
associated to $X$ and $Y$ satisfies $\kappa(t)\geq 0$ when $t\in
[0,\varepsilon)$; if one choice of $\varepsilon$ holds for all pairs
$X,$ $Y,$ then $\Phi_t$ is \emph{locally
nonnegative}.\end{definition}

The notion of an infinitesimally nonnegative endomorphism first
arose in \cite{Tapp}, where Tapp considers the infinitesimally
nonnegative endomorphisms of $\So(4)$ in an effort to classify the
nonnegatively curved metrics on $\SO(4)$.

Given a left-invariant metric $h$ on $G$ with matrix $\Phi$, there
is a unique inverse-linear path $\Phi_t$ with $\Phi_1=\Phi$.  Most
of the inverse-linear paths we are interested in are constructed in
this manner, by specifying the left-invariant metric at time $t=1$.
Such a path is called \emph{nonnegative} if the metric $h_t$ is
nonnegatively curved for $0\leq t\leq 1$.  The following conjecture
of Tapp provides substantial motivation for considering
inverse-linear variations.

\begin{conjecture}[\cite{Tapp}]\label{Tapp Conjecture}
Let $h$ be a nonnegatively curved left-invariant metric on $G$. The
unique inverse-linear path from $h_0$ to $h$ is nonnegative.
\end{conjecture}

This conjecture is particularly interesting since it is not even
known if the space of nonnegatively curved left-invariant metrics on
$G$ is path-connected.  The conjecture is known to be true when $G$
is one of $\SO(3)$ or $\U(2)$ (this follows from the complete
classification of nonnegatively curved left-invariant metrics on
these spaces in \cite{Brown}), and it is true for all known
left-invariant metrics on $\SO(4)$ (\cite{Tapp}). While a proof of
the conjecture eludes us, we do have the following theorem.

\begin{theorem}\label{Main Theorem}
Let $h$ be a nonnegatively curved left-invariant metric on $G$. The
unique inverse-linear path from $h_0$ to $h$ is infinitesimally
nonnegative.
\end{theorem}

The techniques we develop to prove this theorem could potentially be
useful in a proof of Tapp's conjecture.  This theorem is also
interesting in its own right, since it can be used in a
classification of the left-invariant metrics of nonnegative
curvature on $G$.  In light of this theorem, one can first classify
the infinitesimally nonnegative endomorphisms $\Psi$, and then
determine which endomorphisms $(I-\Psi)^{-1}$ correspond to
nonnegatively curved metrics.  Classifying the infinitesimally
nonnegative endomorphisms of $\g$ seems to be an easier question
than classifying the nonnegatively curved left-invariant metrics on
$G$.

We will calculate the Taylor series of $\kappa(t)$ at $0$ to deduce
a formula for $\kappa(t)$ in Section \ref{Taylor Series}, and this
is the essential ingredient to the proof of Theorem \ref{Main
Theorem}. We prove our main theorem in Section \ref{Main Theorem
Section}.  In Section \ref{Examples} we use the Taylor series for
$\kappa(t)$ to classify which subalgebras $\h$ of $\g$ can be
expanded while maintaining nonnegative curvature. In Section
\ref{SO4 Section} we apply our main theorem to Tapp's partial
classification of the infinitesimally nonnegative endomorphisms of
$\So(4)$ to obtain a partial classification of the nonnegatively
curved left-invariant metrics on $\SO(4)$.  We conclude the paper in
Section \ref{Bi-invariant metric shifting} by examining what happens
when the bi-invariant metric on $G$ is changed.

The author would like to thank Kristopher Tapp, Zachary Madden, Nela
Vukmirovic, Angela Doyle, and Min Kim for numerous helpful
discussions and comments on this work.

\section{A formula for $\kappa(t)$}\label{Taylor Series}

Throughout this section, we let $\Phi_t=(I-t\Psi)^{-1}$ be a fixed
inverse-linear variation from the bi-invariant metric $h_0$ on $G$.
We also fix two vectors $X,$ $Y\in \g$.  As described in the
introduction, we obtain a function $\kappa(t)$ determined by
$\Phi_t$, $X$, and $Y$.  The domain of $\kappa$ is the set of all
$t$ such that $\Phi_t$ corresponds to a metric on $G$; this set is
always an open interval of $\R$, determined by the set of
eigenvalues of $\Psi$.  Certain elements of $\g$ will appear
frequently in what follows, so to simplify the exposition we
introduce the Lie algebra elements
\begin{eqnarray*}
A^\Psi &=& [\Psi X,Y]+[X,\Psi Y]\\
B^\Psi &=& [\Psi X,\Psi Y]\\
C^\Psi &=& [\Psi X,Y]+[\Psi Y,X]\\
D^\Psi &=& \Psi^2 [X,Y] -\Psi A^\Psi +B^\Psi.
\end{eqnarray*}
As with $\kappa^\Psi(t)$, we omit the superscripts when no confusion
will arise. Our main tool throughout the article is an explicit
formula for $\kappa(t)$. If $Z_1,$ $Z_2\in \g$, we put $\langle
Z_1,Z_2\rangle = h_0(Z_1,Z_2)$, $|Z_1|^2 = h_0(Z_1,Z_1)$, and
$|Z_1|_{h_t}^2 = h_t(Z_1,Z_1)=\langle \Phi_t Z_1,Z_1\rangle$.

\begin{theorem}\label{kappa formula theorem}
For any $t$ in the domain of $\kappa$, \begin{equation}\label{kappa
formula}\kappa(t) = \alpha + \beta t + \gamma t^2 + \delta t^3 -
\frac 34 t^4\cdot |D|_{h_t}^2,\end{equation} where
\begin{eqnarray*}
\alpha &=& \frac{1}{4} |[X,Y]|^2\\
\beta &=& -\frac{3}{4} \langle \Psi [X,Y],[X,Y]\rangle\\
\gamma &=& -\frac 34
|\Psi[X,Y]|^2+\frac 32 \langle \Psi [X,Y],A\rangle-\frac 12 \langle [X,Y],B\rangle \\
&&-\frac 14 |A|^2 +\frac 14 |C|^2 - \langle [\Psi X,X],[\Psi
Y,Y]\rangle\\
 \delta &=& -\frac 34
 \langle \Psi^3 [X,Y],[X,Y]\rangle+\frac 32 \langle \Psi^2 [X,Y],A
 \rangle - \frac 32 \langle \Psi
 [X,Y],B\rangle   \\&& - \frac 34 \langle \Psi A,A\rangle-\frac 14 \langle \Psi C,
 C\rangle+\langle \Psi[\Psi X,X], [\Psi Y,Y]\rangle+\langle
 A,B\rangle.
\end{eqnarray*}
\end{theorem}

There are two steps to the proof of this theorem.  First we prove
that Equation \ref{kappa formula} holds for all sufficiently small
$t$.  Next we show that each side of the equation is analytic.  This
allows us to invoke the well-known identity theorem: if $f,$
$g\colon I\to \R$ are analytic on an open interval $I$ and $f$ and
$g$ agree on a subinterval of $I$, then $f=g$.  We therefore
conclude that Equation \ref{kappa formula} holds for all $t$. To
accomplish the first step, we calculate the Taylor series of
$\kappa(t)$ at $0$. This calculation will also serve as the
foundation for our analyticity arguments.

\begin{proposition}\label{Power Series}
The Taylor series of $\kappa(t)$ at $0$ is given by $$\kappa(t) =
\alpha+\beta t + \gamma t^2 +\delta t^3 +\Gamma(t),$$ where
$$\Gamma(t) = -\frac 34 \sum_{n=4}^\infty t^n \langle \Psi^{n-4}
D,D\rangle.$$ This formula is valid whenever $|t|< \|\Psi\|^{-1}$,
where $\|\Psi\|=\sup_{|X|=1} |\Psi X|$ is the operator norm of
$\Psi$.
\end{proposition}
\begin{proof}
In \cite{Puttmann}, P\"uttmann shows that the unnormalized sectional
curvature of vectors $Z_1,$ $Z_2\in \g$ with respect to a
left-invariant metric $h$ whose matrix with respect to $h_0$ is
$\Phi$ is given by
\begin{eqnarray}
k_{h}(Z_1,Z_2) &=& \frac 12 \langle [\Phi Z_1,Z_2] + [Z_1,\Phi
Z_2],[Z_1,Z_2]\rangle - \frac 34 |[Z_1,Z_2]|_h^2 \nonumber \\
&&+\frac 14 \langle [Z_1,\Phi Z_2] +[Z_2,\Phi Z_1],\Phi^{-1}
([Z_1,\Phi Z_2] + [Z_2,\Phi Z_1])\rangle \label{Puttmanns formula}\\
&&-\langle [Z_1,\Phi Z_1], \Phi^{-1} [Z_2,\Phi Z_2]\rangle.\nonumber
\end{eqnarray} It follows that
\begin{eqnarray*}
\kappa(t) &=& k_{h_t}(\Phi_t^{-1} X,\Phi_t^{-1} Y)\\
&=& \frac 12 \langle [X,\Phi_t^{-1} Y] +[\Phi_t^{-1}
X,Y],[\Phi_t^{-1}X,\Phi_t^{-1}Y]\rangle\\
&&-\frac 34 \langle \Phi_t [\Phi_t^{-1} X,\Phi_t^{-1}Y],[\Phi_t^{-1}
X,\Phi_t^{-1} Y]\rangle\\
&&+\frac 14 \langle [\Phi_t^{-1} X,Y] + [\Phi_t^{-1}
Y,X],\Phi_t^{-1} ([\Phi_t^{-1} X,Y]+[\Phi_t^{-1} Y,X])\rangle\\
&& - \langle [\Phi_t^{-1} X,X],\Phi_t^{-1} [ \Phi_t^{-1}
Y,Y]\rangle\\
&=& I_1 - I_2 + I_3 - I_4.
\end{eqnarray*}
Using the expression $\Phi_t^{-1} = I-t\Psi$, we can easily simplify
$I_1$, $I_3$, and $I_4$.  We find
\begin{eqnarray*}
I_1 &=& |[X,Y]|^2 - \frac {3t}2 \langle [X,Y],A\rangle + t^2(\langle
[X,Y],B\rangle +\frac 12 |A|^2 \rangle - \frac {t^3}2 \langle
A,B\rangle\\
I_3 &=& \frac {t^2}4 |C|^2 - \frac {t^3}4\langle C,\Psi C\rangle\\
I_4 &=& t^2 \langle [\Psi X,X],[\Psi Y,Y]\rangle - t^3 \langle [\Psi
X,X],\Psi[\Psi Y,Y]\rangle.
\end{eqnarray*}
To calculate $I_2$, notice that if $|t|< \|\Psi\|^{-1}$, then
$$\Phi_t = \sum_{n=0}^\infty t^n \Psi^n,$$ where the convergence is
in the space of endomorphisms of $\g$ with the operator norm.  Using
this expression we calculate
\begin{eqnarray*}
\frac 43 I_2 &=&  \langle \Phi_t([X,Y]-t A +t^2 B),[X,Y]-t A +t^2
B\rangle\\
&=&  \sum_{n=0}^\infty t^n \langle \Psi^n [X,Y] - t \Psi^n A
+ t^2 \Psi^n B,[X,Y]-t A +t^2 B\rangle\\
&=&  \sum_{n=0}^\infty t^n \left(\langle \Psi^n[X,Y],[X,Y]\rangle-2t
\langle \Psi^n [X,Y],A\rangle \right.\\&&\left.+t^2(\langle \Psi^n
A,A\rangle+2\langle \Psi^n [X,Y],B\rangle)-2t^3 \langle \Psi^n
A,B\rangle+t^4 \langle \Psi^n B,B\rangle \right)
\end{eqnarray*}

\begin{eqnarray*}&=& |[X,Y]|^2 +t(\langle \Psi [X,Y],[X,Y]\rangle  - 2 \langle
[X,Y],A\rangle)\\&&+t^2(\langle \Psi^2 [X,Y],[X,Y]\rangle-2\langle
\Psi [X,Y],A\rangle+|A|^2+2\langle [X,Y],B\rangle)\\
&&+t^3(\langle \Psi^3[X,Y],[X,Y]\rangle - 2 \langle
\Psi^2[X,Y],A\rangle +\langle \Psi A,A\rangle \\&&\quad+2\langle
\Psi [X,Y],B\rangle-2\langle A,B\rangle)
\\&&+\sum_{n=4}^\infty t^n \langle \Psi^{n-4} D,D\rangle.\end{eqnarray*}

Combining the different terms proves the result.
\end{proof}

Notice that the power series of $\kappa(t)$ would have been much
messier if we were considering the unnormalized sectional curvature
of $X$ and $Y$ with respect to $h_t$ instead of the unnormalized
sectional curvature of $\Phi_t^{-1} X$ and $\Phi_t^{-1} Y$.

When $|t|<\|\Psi\|^{-1}$, observe that $$\Gamma(t) = -\frac 34
\sum_{n=4}^\infty t^n\langle \Psi^{n-4} D,D\rangle  = -\frac 34 t^4
\langle \Phi_t D,D\rangle = -\frac 34 t^4 \cdot |D|_{h_t}^2.$$ This
proves that Equation \ref{kappa formula} holds for small $t$.
Therefore to complete the proof of Theorem \ref{kappa formula
theorem}, all we must do is prove that $\kappa(t)$ and $|D|_{h_t}^2$
are analytic.

\begin{lemma}\label{analyticity}
The function $\kappa(t)$ is analytic on its domain of definition.
\end{lemma}
\begin{proof}
Assume that $t_0$ is such that $\Phi_{t_0}$ corresponds to a metric
on $G$. We show that $\kappa$ is locally a power series at $t_0$.
Recalling P\"utmann's Formula \ref{Puttmanns formula}, it is clear
that we must only prove that
$$|[\Phi_t^{-1} X,\Phi_t^{-1} Y]|^2_{h_t}$$ can be expressed as a power
series near $t_0$.  Since $\Psi$ is $h_0$-self-adjoint, it can be
diagonalized; say $\Psi=\diag(a_1,\ldots,a_d)$.  We then have
\begin{eqnarray}\label{Phi t equation} \Phi_t&=&\diag\left(\frac
1{1-a_1t},\ldots,\frac 1{1-a_dt}\right)\nonumber\\
&=& \diag\left(\frac{1}{1-a_it_0} \sum_{n=0}^\infty
\left(\frac{a_i}{1-a_i t_0}\right)^n
(t-t_0)^n\right)_{i=1}^d\\
&=& \Phi_{t_0} \sum_{n=0}^\infty \Phi_{t_0}^n \Psi^n
(t-t_0)^n\nonumber,\end{eqnarray} with convergence whenever
$|t-t_0|$ is sufficiently small.  We can use this expression for
$\Phi_t$ together with the identity $\Phi_t^{-1} = I-t_0 \Psi
-(t-t_0)\Psi$ to expand $|[\Phi_t^{-1} X,\Phi_t^{-1} Y]|_{h_t}^2$ as
a power series as in the proof of Proposition \ref{Power Series}.
\end{proof}

The analyticity of $|D|_{h_t}^2$ also follows from Equation \ref{Phi
t equation}, completing the proof of Theorem \ref{kappa formula
theorem}.

\section{The main theorem}\label{Main Theorem Section}

We now use the formula for $\kappa(t)$ from the previous section to
prove our main theorem.  First, we use the power series of
$\kappa(t)$ to rephrase what it means for $\Psi$ to be
infinitesimally nonnegative.

\begin{proposition}\label{Inf Nonneg Equivalence}
An $h_0$-self-adjoint endomorphism $\Psi\colon \g\to\g$ is
infinitesimally nonnegative if and only if whenever $X,$  $Y\in \g$
commute we have either
\begin{enumerate}
\item $\kappa'''(0) \geq 0$ or
\item $\kappa'''(0)=0$ and $D=0$.
\end{enumerate}
In the second case, $\kappa$ is identically zero.
\end{proposition}

Notice that $\kappa'''(0)=6\delta$, so that $\kappa'''(0)$ is a
relatively simple sum of inner products of elements of $\g$.  This
equivalent definition of infinitesimally nonnegative is often more
useful than our original definition, especially in applications.

\begin{proof}
It is clear that $\Psi$ is infinitesimally nonnegative if and only
if either all coefficients of the Taylor series of $\kappa(t)$ at
$0$ are zero or the first nonzero coefficient is positive.  When $X$
and $Y$ do not commute, we have $\alpha>0$.

Suppose instead that $X,$ $Y\in \g$ commute. Notice that
$\alpha=\beta=\gamma=0$; in particular $\gamma=0$ follows from the
observation that $$-\frac 14 |A|^2 + \frac 14 |C|^2 -\langle [\Psi
X,X],[\Psi Y,Y]\rangle =0,$$ a consequence of the Jacobi identity
and the bi-invariance of $h_0$.  If $\delta>0$, then the first
nonzero coefficient of the power series is positive. The coefficient
of $t^4$ is $-(3/4)\cdot |D|^2$, so if $\delta=0$ we must have $D=0$
in order for $\Psi$ to be infinitesimally nonnegative.  When $D=0$,
the coefficients of $t^n$ for $n>4$ are also $0$, so these
conditions are also sufficient to ensure that $\Psi$ is
infinitesimally nonnegative.
\end{proof}

With this equivalent definition of infinitesimally nonnegative, we
are prepared to prove our main theorem.

\begin{mainTheorem}
Let $h$ be a nonnegatively curved left-invariant metric on $G$.  The
unique inverse-linear path from $h_0$ to $h$ is infinitesimally
nonnegative.
\end{mainTheorem}
\begin{proof}
Let $\Phi_t$ be the unique inverse-linear path from $h_0$ to $h$,
and fix commuting $X,$ $Y\in \g$. By Theorem \ref{kappa formula
theorem},
$$\kappa(t) = \frac 16 \kappa'''(0)t^3 -\frac 34 t^4\cdot
|D|_{h_t}^2.$$ As $h$ is nonnegatively curved, $$0\leq \kappa(1) =
\frac 16 \kappa'''(0) -\frac 34 \cdot |D|_{h}^2.$$ Therefore
$\kappa'''(0)\geq 0$, and if $\kappa'''(0)=0$ then $D=0$.  By
Proposition \ref{Inf Nonneg Equivalence}, $\Psi$ is infinitesimally
nonnegative.
\end{proof}

As an example of how Theorem \ref{Main Theorem} can be applied to
transform infinitesimal results into global ones, consider the
following result from \cite{Tapp}.

\begin{lemma}[{\cite[Lemma 2.6]{Tapp}}]\label{infinitesimal rigidity}
Assume that $\Psi$ is infinitesimally nonnegative.  Let $\p_0$ be
the eigenspace of $\Psi$ corresponding to the smallest eigenvalue.
If $X\in \p_0$, $Y\in \g$ and $[X,Y]=0$, then $[X,\Psi Y]\in \p_0$.
\end{lemma}

This lemma was the main tool used by Tapp to prove rigidity
statements about the family of infinitesimally nonnegative
endomorphisms of $\So(4)$.  Applying Theorem \ref{Main Theorem}, we
derive the following global result.

\begin{lemma}\label{global rigidity}
Assume that $\Phi$ is the matrix of a nonnegatively curved metric.
Let $\p_0$ be the eigenspace of $\Phi$ corresponding to the smallest
eigenvalue of $\Psi$.  If $X\in \p_0$, $Y\in \g$ and $[X,Y]=0$, then
$[X,\Phi^{-1} Y]\in \p_0$.
\end{lemma}
\begin{proof}
This follows immediately from Lemma \ref{infinitesimal rigidity}
since $\Psi=I-\Phi^{-1}$ is infinitesimally nonnegative and $\p_0$
is the eigenspace of $\Psi$ corresponding to the smallest
eigenvalue.
\end{proof}

We note that this result can also be derived directly from
Puttmann's Formula \ref{Puttmanns formula}.

\section{Enlarging subalgebras}\label{Examples}

Perhaps the simplest type of inverse-linear variation is one which
gradually scales vectors in a subalgebra $\h$ of $\g$.  In this
section, we let $H\subset G$ be a Lie subgroup of the Lie group $G$
with Lie algebra $\h\subset \g$.  It is known that shrinking vectors
in $\h$ yields nonnegatively curved metrics (see, for example,
\cite{Eschenburg}). On the other hand, expanding vectors in $\h$
does not always produce nonnegatively curved metrics, but does when
$\h$ is abelian \cite{Grove}. In this section, we use the power
series formula from the previous section to determine when the
subalgebra $\h$ can be expanded while maintaining nonnegative
curvature.

For $Z\in \g$, denote by $Z^\h$ and $Z^\p$ the projections of $Z$
onto $\h$ and its $h_0$-orthogonal complement $\p$.  When we discuss
scaling $\h$ by a factor $\lambda>0$, we are referring to the metric
$$h(X,Y) = h_0(\lambda X^\h + X^\p,Y),$$ so that the square of the norm of
a vector in $h$ is scaled by a factor of $\lambda$.  Let
$\Psi(Z)=Z^\h$, so $\Phi_t = (I-t\Psi)^{-1}$ is the inverse-linear
variation which gradually expands vectors in $\h$. If $\h$ is
abelian, it is easy to use the formulas for the coefficients of the
power series of $\kappa(t)$ in tandem with the analyticity of
$\kappa$ to prove
\begin{equation}\label{Abelian blowup curv eq}\kappa(t) = \frac 14 |[X,Y]|^2 - \frac 34 |[X,Y]^\h|^2
\cdot \frac t{1-t} \qquad (-\infty < t < 1).\end{equation} From this
formula we can show that enlarging $\h$ by a factor of up to $4/3$
always preserves nonnegative curvature, a result which first
appeared in \cite{Grove}.  In fact, the particularly nice form of
$\kappa(t)$ allows us to prove a stronger statement.

\begin{theorem}\label{Abelian Blowup}
Scaling the abelian subalgebra $\h\subset \g$ preserves nonnegative
curvature if and only if no vector in $[\g,\g]$ is expanded by more
than $4/3$.
\end{theorem}
\begin{proof}
By Equation \ref{Abelian blowup curv eq}, the metric $h_t$ is
nonnegatively curved if and only if
\begin{equation}\label{blowup inequality} |Z^\h|^2 \cdot \frac t {1-t} \leq \frac 13|Z|^2
 \end{equation} holds for all $Z\in [\g,\g]$. As
 $$|Z|_{h_t}^2 = \langle \Phi_t Z,Z\rangle = \langle Z + \frac
 t{1-t} Z^\h,Z\rangle = |Z|^2 + |Z^\h|^2 \cdot \frac t{1-t},$$ we
 find that Inequality \ref{blowup inequality} is equivalent to
 requiring that $|Z|_{h_t}^2 \leq (4/3)\cdot |Z|^2$ holds for all
 $Z\in [\g,\g]$.
\end{proof}

If $[\g,\g]\cap \h\neq 0$, this theorem says that $\h$ can be by a
factor up to $4/3$. At the other extreme, if $[\g,\g]\perp \h$ then
we find that $\h$ can be expanded up by an arbitrary amount.  This
was already known, since if $\h$ is orthogonal to $[\g,\g]$ then
$\h$ is contained in the center of $\g$.  This rescaling then stays
within the family of bi-invariant metrics on $\g$.

When $\h$ is not abelian, things are not quite so simple.  In this
case the power series simplifies to
$$\kappa(t) = \frac 14 |[X,Y]|^2 -\frac 34 |[X,Y]^\h|^2 t +\frac 34
|B|^2 t^2 - \frac 14 |B|^2 t^3 -\frac 34 |[X^\p,Y^\p]^\h|^2 \cdot
\frac {t^2}{1-t}.$$ We can use this formula to classify exactly
which subalgebras of $\g$ can be enlarged a small amount while
maintaining nonnegative curvature.

\begin{theorem}\label{Arbitrary Blowup}
Expanding the subalgebra $\h\subset \g$ by a small amount preserves
nonnegative curvature if and only if there exists a constant $c$
such that $|[X^\h,Y^\h]|\leq c\cdot |[X,Y]|$ holds for all $X,$
$Y\in \g$.
\end{theorem}

We omit the lengthy but easy proof for the reason that we do not
know if there are any interesting examples of subalgebras for which
the latter condition holds.  It clearly holds when $\h$ is either
abelian or an ideal of $\g$ (or the sum of an ideal and a disjoint
abelian subalgebra), but it is already known that such subalgebras
can be enlarged while maintaining nonnegative curvature.

\section{An application to $\SO(4)$}\label{SO4 Section}

In this section we apply our main theorem to Tapp's partial
classification of the infinitesimally nonnegative endomorphisms of
$\So(4)$ to give a partial classification of the nonnegatively
curved left-invariant metrics on $\SO(4)$.  Let $G=\SO(4)$, so
$\g=\So(4)=\So(3)\oplus\So(3)=\g_1\oplus\g_2$.  When we give
$\So(3)$ the bi-invariant metric $$\langle Z_1,Z_2\rangle = \frac
12\cdot \trace(Z_1 Z_2^T),$$ $\g$ inherits the bi-invariant product
metric $h_0$ from its factors.  We call a vector in $\g$
\emph{singular} if it is in either $\g_1$ or $\g_2$.

The known examples of left-invariant metrics of nonnegative
curvature are discussed thoroughly in \cite[Section 3]{Tapp}.  These
metrics can be grouped into three categories:  product metrics,
metrics which come from a torus action, and metrics which come from
an $S^3$-action.  Our main result in this section is the following
theorem.

\begin{theorem}\label{Singular eigenvector classification}
Let $\Phi$ be the matrix of a nonnegatively curved left-invariant
metric $h$.  If $\Phi$ has a singular eigenvector, then either $h$
is a product metric or $h$ comes from a torus action.  In either
case, $h$ is a known example of a metric of nonnegative curvature.
\end{theorem}

Since the nonnegatively curved left-invariant metrics on $\SO(3)$
have been classified, the product metrics are well-understood.

The metrics arising from torus actions are only a little more
complicated than the product metrics.  Let $\{A_1,A_2,A_3\}$ and
$\{B_1,B_2,B_3\}$ be $h_0$-orthonormal bases of $\g_1$ and $\g_2$,
respectively.  After scaling $\g_1$ and $\g_2$ by factors $c$ and
$d$, respectively, then enlarging the abelian subalgebra
$\tau=\spn\{A_3,B_1\}$ by $4/3$, then further altering the metric on
$\tau$ via a canonical $T^2$-action on $G$ (using the method of
Cheeger \cite{Cheeger}), one obtains a nonnegatively curved metric
$h$ with matrix
\begin{equation}\label{Torus Form}\Phi = \left(\begin{array}{cccccc} c & 0 & 0 & 0 & 0 & 0\\
0 & c & 0 & 0 & 0&0\\
0 & 0 &a_1 & a_3 & 0 & 0\\
0 & 0 &a_3 & a_2 & 0 & 0\\
0 & 0 & 0 & 0 & d &0\\
0 &0 &0 & 0 & 0 & d
\end{array}\right)\end{equation}
with respect to the basis $\{A_1,A_2,A_3,B_1,B_2,B_3\}$.  The only
restriction on $\Phi$, coming from the fact that the final
alteration only shrinks vectors (see \cite[Section 3.2]{Tapp}), is
that the norm on $\tau$ determined by the matrix
$$\left(\begin{array}{cc} a_1 & a_3 \\ a_3 & a_2\end{array}\right)$$
is bounded above by the norm determined by $$\left(\begin{array}{cc}
\frac 43 \cdot c & 0 \\ 0 & \frac 43 \cdot d\end{array}\right).$$

The most difficult part of the proof of Theorem \ref{Singular
eigenvector classification} was already completed in \cite[Theorem
4.1]{Tapp}, which is an infinitesimal version of the theorem. We
call an endomorphism $\Psi$ of $\g$ a \emph{product} if $\Psi
(\g_1)\subset \g_1$ and $\Psi(\g_2)\subset \g_2$.

\begin{theorem}[{\cite[Theorem 4.1]{Tapp}}]\label{inf singular eigen classification}
Let $\Psi$ be an infinitesimally nonnegative endomorphism with a
singular eigenvector.  Then either $\Psi$ is a product or $\Psi$ is
of Form \ref{Torus Form}.
\end{theorem}

\begin{proof}[Proof of Theorem \ref{Singular eigenvector
classification}] Since $h$ is nonnegatively curved,
$\Psi=I-\Phi^{-1}$ is infinitesimally nonnegative.  As $\Phi$ has a
singular eigenvector, so does $\Psi$.  According to Theorem \ref{inf
singular eigen classification}, either $\Psi$ is a product or $\Psi$
can be written in Form \ref{Torus Form}.  Clearly if $\Psi$ is a
product then $\Phi$ is a product, which means $h$ is a product
metric.  If instead $\Psi$ has Form \ref{Torus Form}, then so does
$\Phi$.

Assume $\Phi$ has Form \ref{Torus Form}; we must prove that $\Phi$
comes from a torus action. Permuting some basis vectors if
necessary, we may assume that $A_1$, $A_2$, $A_3$ and $B_1$, $B_2$,
$B_3$ behave like the quaternions $\B i$, $\B j$, $\B k$ with
respect to their Lie bracket structure. Denote by $\tilde h$ the
metric on $\tau$ corresponding to the matrix $$
\left(\begin{array}{cc} \frac 43 \cdot c & 0 \\ 0 & \frac 43 \cdot
d\end{array}\right).$$ We must prove that
\begin{equation*}|\alpha A_3+\beta B_1|_h^2
\leq |\alpha A_3 +\beta B_1|_{\tilde h}^2\end{equation*} holds for
all $\alpha,\beta \in \R$.

Consider the unnormalized sectional curvature of the vectors $\alpha
A_1+\beta B_2$ and $A_2+B_3$ with respect to $h$. We have
\begin{eqnarray*} [\Phi (\alpha A_1+\beta B_2),A_2+B_3] &=& \alpha c
A_3+\beta d B_1\\{} [\alpha A_1+\beta B_2,\Phi(A_2+B_3)] &=& \alpha
c A_3 + \beta d B_1\\{} [\alpha A_1 + \beta B_2,A_2+B_3] &=& \alpha
A_3 +\beta B_1,
\end{eqnarray*}
and therefore by P\"uttmann's Formula \ref{Puttmanns formula}
\begin{eqnarray*}
k_h(\alpha A_1+\beta B_2,A_2+B_3) &=& \langle \alpha c A_3 + \beta d
B_1,\alpha A_3+\beta B_1\rangle-\frac 34 |\alpha A_3 + \beta B_1|_h^2\\
&=& \frac 34 (|\alpha A_3+\beta B_1|_{\tilde h}^2 - |\alpha A_3
+\beta B_1|_{h}^2).
\end{eqnarray*}
Since $h$ is nonnegatively curved, this proves the required
inequality.
\end{proof}

We have therefore completed the classification of nonnegatively
curved metrics with singular eigenvectors.  The rest of the
classification is still a difficult problem.  Among all the
remaining left-invariant metrics, we must still locate the metrics
which come from $S^3$-actions.  Noticing that there are three
$\Phi$-invariant $2$-dimensional abelian subalgebras of $\g$
whenever $\Phi$ comes from an $S^3$-action (\cite[Section
3.3]{Tapp}), Tapp looked for such subalgebras in the infinitesimal
version of the problem.

\begin{theorem}[{\cite[Theorem 4.4]{Tapp}}]\label{Tapp SO(4) Theorem}
If $\Psi$ is an infinitesimally nonnegative endomorphism of $\g$,
then $\g$ has a $\Psi$-invariant $2$-dimensional abelian subalgebra.
\end{theorem}

Applying our main theorem, we immediately have the following result.

\begin{theorem}\label{2-dim subspace}
If $h$ is nonnegatively curved, then $\g$ has a $\Phi$-invariant
$2$-dimensional abelian subalgebra.
\end{theorem}

We can actually say slightly more, by mimicking the argument of
 \cite[Theorem 4.2]{Tapp}.

\begin{theorem}
There are orthonormal bases $\{A_1,A_2,A_3\}$ and $\{B_1,B_2,B_3\}$
of the two factors of $\g=\g_1\oplus \g_2$ such that with respect to
the basis $\{A_1,B_1,A_2,B_2,A_3,B_3\}$, $\Phi$ has the form
$$
\Phi= \left(\begin{array}{cccccc} a_1 & a_3 & 0 & 0 & 0 & 0\\ a_3 &
a_2 & 0 & 0 & 0 & 0 \\ 0 & 0 & b_1 & b_3 & \lambda & 0\\ 0 & 0 & b_3
& b_2 & 0 & \mu\\ 0 & 0 & \lambda & 0 & c_1 & c_3 \\ 0 & 0 & 0 & \mu
& c_3 & c_2 \end{array}\right).$$
\end{theorem}
\begin{proof}
In \cite[Theorem 4.2]{Tapp}, Tapp proves that $\Psi$ can be written
in this form. The only properties of $\Psi$ which are used are that
there is a $2$-dimensional $\Psi$-invariant abelian subalgebra and
that $\Psi$ is self-adjoint.  Since Theorem \ref{2-dim subspace}
implies that $\Phi$ has these properties as well, the proof carries
over to $\Phi$ immediately.
\end{proof}

If it is the case that $\lambda=\mu=0$, then there are three
$\Phi$-invariant $2$-dimensional abelian subalgebras.

\section{Changing the bi-invariant metric}\label{Bi-invariant metric shifting}

Suppose that the unique inverse-linear path from the bi-invariant
metric $h_0$ to the left-invariant metric $h$ is nonnegative.  If
$h_1$ is another bi-invariant metric, must the unique inverse-linear
path from $h_1$ to $h$ also be nonnegative?  If Tapp's Conjecture
\ref{Tapp Conjecture} is true, then the answer must be yes.  A proof
of this statement could also be a logical step towards proving the
conjecture, as it would imply that Conjecture \ref{Tapp Conjecture}
is equivalent to the seemingly weaker statement that if $h$ is
nonnegatively curved then there is some bi-invariant metric $h_0$
such that the unique inverse-linear path from $h_0$ to $h$ is
nonnegative. In this section, we address the case where $h_1$ is a
scalar multiple of $h_0$.

\begin{theorem}\label{transferable properties}
If $\lambda>0$ and $h$ is any left-invariant metric, then
\begin{enumerate}
\item the unique inverse-linear path from $h_0$ to $h$ is
nonnegative if and only if the unique inverse-linear path from
$\lambda h_0$ to $h$ is nonnegative.

\item The unique inverse-linear path from $h_0$ to $h$ is locally
nonnegative if and only if the unique inverse-linear path from
$\lambda h_0$ to $h$ is locally nonnegative.

\item The unique inverse-linear
path from $h_0$ to $h$ is infinitesimally nonnegative if and only if
the unique inverse-linear path from $\lambda h_0$ to $h$ is
infinitesimally nonnegative.
\end{enumerate}
\end{theorem}

Let $\Phi$ be the matrix of the left-invariant metric $h$ with
respect to $h_0$.  Let $\lambda>0$, and let $\Theta$ be the matrix
of $h$ with respect to the bi-invariant metric $\lambda h_0$.  When
we let $\Psi = I-\Phi^{-1}$ and $\Upsilon = I-\Theta^{-1}$, we see
that $\Phi_t = (I-t\Psi)^{-1}$ is the unique inverse-linear path
from $h_0$ to $h$ and $\Theta_t=(I-t\Upsilon)^{-1}$ is the unique
inverse-linear path from $\lambda h_0$ to $h$.  As
$$h_0(\Phi X,Y)=h(X,Y) = \lambda h_0(\Theta X,Y),$$ we find that
$\Theta = \lambda^{-1} \Phi$. Theorem \ref{transferable properties}
will follow immediately from our next proposition.

\begin{proposition}\label{Curvature relation}
With the notation of the previous paragraph,
\begin{equation}\label{curv rel eq}\kappa^{\Upsilon}(t) = \lambda (1-(1-\lambda)t)^3
\cdot \kappa^\Psi\left(\frac{ \lambda t}{1-(1-\lambda)t}\right)
\qquad (0\leq t \leq 1).\end{equation}
\end{proposition}
\begin{proof}
By analyticity, it is enough to show that the above equality holds
for all sufficiently small $t$.  We therefore assume that $t$ is
small enough that the various power series which appear in the
following proof are convergent. Since $\Theta=\lambda ^{-1}\Phi$, we
find
\begin{equation*} \Upsilon = I-\lambda \Phi^{-1} = I-\lambda (I-\Psi) =
(1-\lambda )I + \lambda \Psi.\end{equation*} A very long, but
entirely straightforward calculation using this equality now shows
\begin{eqnarray}
\alpha^\Upsilon &=& \lambda  \alpha^\Psi\nonumber\\
\beta^\Upsilon &=& -3(1-\lambda )\lambda \alpha^\Psi +\lambda ^2 \beta^\Psi\nonumber\\
\gamma^\Upsilon &=& 3(1-\lambda )^2 \lambda  \alpha^\Psi - 2(1-\lambda )\lambda ^2 \beta^\Psi + \lambda ^3 \gamma^\Psi\nonumber\\
\delta^\Upsilon &=& -(1-\lambda )^3 \lambda  \alpha^\Psi +
(1-\lambda )^2 \lambda ^2 \beta^\Psi - (1-\lambda )\lambda ^3
\gamma^\Psi +\lambda ^4 \delta^\Psi,\label{delta equality}
\end{eqnarray}
from which we conclude that
\begin{eqnarray*}\label{another curv relation eq}
\lambda ^{-1} \kappa^\Upsilon(t) &=& \alpha^\Psi(1-(1-\lambda )t)^3
+ \beta^\Psi t \lambda  (1-(1-\lambda )t)^2 \\&&+ \gamma^\Psi t^2
\lambda ^2 (1-(1-\lambda )t) + \delta^\Psi t^3 \lambda ^3 +\lambda
^{-1}\Gamma^\Upsilon(t).\end{eqnarray*}
 Dividing both sides of this
equation by $(1-(1-\lambda )t)^3$, we see
$$\frac{\kappa^\Upsilon(t)}{\lambda (1-(1-\lambda )t)^3} =
\kappa^\Psi\left(\frac {\lambda  t}{1-(1-\lambda )t}\right) -
\Gamma^\Psi\left(\frac {\lambda t}{1-(1-\lambda
)t}\right)+\frac{\Gamma^\Upsilon(t)}{\lambda (1-(1-\lambda )t)^3},$$
and we must only prove that the last two terms above cancel. Now
$D^\Upsilon = \lambda^2 D^\Psi$, so applying the binomial theorem to
$\Upsilon^{n-4}=((1-\lambda)I+\lambda \Psi)^{n-4}$ we have
\begin{eqnarray*}
\Gamma^\Upsilon(t)&=& \lambda ^5\sum_{n=4}^\infty t^n \langle \Upsilon^{n-4}D^\Psi,D^\Psi \rangle\\
\\&=& \lambda ^5 \sum_{n=4}^\infty \sum_{k=0}^{n-4} t^n {n-4 \choose k} (1-\lambda )^{n-4-k} \lambda ^k \langle \Psi^k D^\Psi,D^\Psi\rangle.\\
\end{eqnarray*}
This double sum is absolutely convergent for all sufficiently small
$t$. Thus we may interchange the order of summation, and we find
\begin{eqnarray*} \Gamma^{\Upsilon}(t) &=& \lambda  \sum_{k=0}^\infty t^{k+4}\lambda ^{k+4}\langle \Psi^k D^\Psi ,D^\Psi\rangle \sum_{n=0}^\infty t^{n} {n+k \choose k} (1-\lambda )^{n}\\
&=& \lambda  \sum_{k=0}^\infty t^{k+4}\lambda ^{k+4}\langle \Psi^k D^\Psi,D^\Psi\rangle \left(\sum_{n=0}^\infty t^n(1-\lambda )^n\right)^{k+1}\\
&=& \lambda  \sum_{k=0}^\infty t^{k+4}\lambda ^{k+4}\langle \Psi^k
D^\Psi,D^\Psi\rangle \cdot \frac{1}{(1-(1-\lambda )t)^{k+1}}.
\end{eqnarray*}
From this we conclude
$$\frac{\Gamma^\Upsilon(t)}{\lambda (1-(1-\lambda )t)^3} =
\Gamma^\Psi\left(\frac{\lambda  t}{1-(1-\lambda )t}\right)$$ holds
for all sufficiently small $t$, as was to be shown.
\end{proof}

\begin{proof}[Proof of Theorem \ref{transferable properties}]
This follows from Proposition \ref{Curvature relation} since
$$t\mapsto \frac{\lambda t}{1-(1-\lambda)t}$$ is an
endpoint-fixing self-homeomorphism of the unit interval.
\end{proof}

While our proof of Theorem \ref{transferable properties} is nice in
that it proves all three parts of the theorem at once, it also hides
the original intuition behind the argument.  Notice that when $X$
and $Y$ commute, Equation \ref{delta equality} implies that
$\delta^\Upsilon = \lambda^4 \delta^\Psi$. Since also
$D^\Upsilon=\lambda^2 D^\Psi$, we may conclude that $\Psi$ is
infinitesimally nonnegative if and only if $\Upsilon$ is
infinitesimally nonnegative by appealing to Proposition \ref{Inf
Nonneg Equivalence}.  This proves the third part of Theorem
\ref{transferable properties}.  We originally found this relation
between $\delta^\Upsilon$ and $\delta^\Psi$, which suggested that we
look for a relation between $\kappa^\Upsilon$ and $\kappa^\Psi$.
This culminated in deriving Equation \ref{curv rel eq} for all
sufficiently small $t$, which is enough to prove the second part of
the theorem. Finally we noticed that $\kappa(t)$ is analytic, which
proves the first part of the theorem.

\bibliographystyle{amsplain}

\begin{thebibliography}{10}

\bibitem{Brown} Brown, Finck, Spencer, Tapp, Wu, \emph{Invariant metrics
with nonnegative curvature on compact Lie groups}, Canadian Math.
Bull., to appear.

\bibitem{Cheeger} J. Cheeger, \emph{Some examples of manifolds of
nonnegative curvature}, J. Differential Geom. \textbf{8} (1972),
623-628.

\bibitem{Eschenburg} J.-H. Eschenburg, \emph{Inhomogeneous spaces of
positive curvature}, Differential Geom. Appl. \textbf{2} (1992).

\bibitem{Grove} K. Grove, W. Ziller, \emph{Curvature and symmetry of
Milnor spheres}, Annals of Math. \textbf{152} (2000) 331-367.

\bibitem{Puttmann} T. P\"uttmann, \emph{Optimal pinching constants of odd
dimensional homogeneous spaces}, Ph. D. thesis, Ruhr-Universit\"at,
Germany, 1991.

\bibitem{Tapp} K. Tapp, \emph{Invariant metrics on $\SO(4)$ with
nonnegative curvature}, Preprint.

\end{thebibliography}

\end{document}